\documentclass{amsart}

\usepackage{amsmath, amssymb, graphicx, psfrag}
\usepackage[mathscr]{eucal}

\newtheorem{thm}{Theorem}[section]
\newtheorem{lem}[thm]{Lemma}
\newtheorem{prop}[thm]{Proposition}

\theoremstyle{definition}

\newtheorem{rem}[thm]{Remark}

\newcommand{\norm}[1]{\parallel\!#1\!\parallel}

\newcommand{\Z}{\ensuremath{\mathbf{Z}}}
\newcommand{\R}{\ensuremath{\mathbf{R}}}

\newcommand{\Hyp}{\ensuremath{\mathbf{H}}}

\def\co{\colon\thinspace}

\begin{document}

\title[Geometrically infinite surfaces]{Geometrically infinite surfaces in 3--manifolds\\ with hyperbolic fundamental group}

\author[J.~Barnard]{Josh Barnard}
\address{Dept.\ of Mathematics\\
        University of Oklahoma\\
	Norman, OK 73019}
\email{jbarnard@math.ou.edu}

\date{\today}

\begin{abstract}
We prove a partial generalization of Bonahon's tameness result to surfaces inside irreducible 3--manifolds with hyperbolic fundamental group.

Bonahon's result states that geometrically infinite ends of freely indecomposable hyperbolic 3--manifolds are simply degenerate. It is easy to see that a geometrically infinite end gives rise to a sequence of curves on the corresponding surface whose geodesic representatives are not contained in any compact set. The main step in his proof is showing that one may assume that these curves are simple on the surface.

In this paper, we generalize the main step of Bonahon's proof, showing that a geometrically infinite end gives rise to a sequence of simple surface curves whose geodesic representatives are not contained in any compact set.
\end{abstract}

\maketitle

\tableofcontents

\section{Introduction}\label{intro}

In~\cite{bonahon:ends} Bonahon shows that geometrically infinite ends of freely indecomposable hyperbolic 3--manifolds are simply degenerate. Suppose $N$ is a closed irreducible hyperbolic 3--manifold containing a $\pi_1$-injective closed surface $S$. (Bonahon's theorem covers a broader class of 3--manifolds, but this is the relevant set-up for our purposes.) Let $M$ be the cover of $N$ corresponding to the subgroup $\pi_1(S)\leq\pi_1(N)$. If $S$ is geometrically infinite in $N$, then there is a sequence of essential curves $\{\gamma_i\}$ on $S$ whose corresponding geodesic representatives $\{\gamma_i^*\}$ in $M$ are not contained in any compact set in $M$ (i.e., they \emph{reach out} an end of $M$). The most important and most technical part of Bonahon's proof is in showing that we may assume that each $\gamma_i$ is simple on $S$. The next step is to find a sequence of pleated surfaces $S_i$ ``hanging off'' these $\gamma_i^*$ which, having uniformly bounded diameter in $M$, must exit an end of $M$. In particular, the simple curves $\gamma_i$ have geodesic representatives $\gamma_i^*$ in $M$ which exit an end of $M$. This is what it means for an end of $M$ to be simply degenerate. Moreover, with this more restrictive set-up, it is easy then to deduce that $S$ is a virtual fiber in $N$. In particular, $M$ is topologically tame (so that, in this case, $M$ is homeomorphic to $S\times\R$).

In this paper we generalize the main step (of the restricted version described above) of Bonahon's proof to the case that $N$ is not known to be hyperbolic, but rather is known only to have Gromov-hyperbolic fundamental group.

There are two important comments to be made. The first is that this result is an immediate consequence of geometrization, as any geometric $N$ satisfying the given hypotheses must support a complete hyperbolic structure, so that Bonahon's theorem applies. On the other hand, our version (cf.~Prop~\ref{prop:ending}) of the main technical proposition of Bonahon's paper~\cite[Prop $5.1$]{bonahon:ends} applies not only to 3--manifolds, but to closed surface subgroups of any torsion-free hyperbolic group $G$. The assumption that $G=\pi_1(N)$ is used only in the intersection lemma~\ref{lem:intnum}. In fact, almost all the novelty of this paper can be found in the proof of Proposition~\ref{prop:ending}.

Second, we note that the pleated surface portion of Bonahon's proof is relatively simple, especially when compared to the subtle arguments contained elsewhere in his paper. In contrast, generalizing this argument to obtain a corresponding statement about hyperbolic groups seems less than straightforward. We will say more about this issue in the next subsection, once the appropriate notation is introduced.

\subsection{Definitions and Notation}\label{basic}

We introduce some notation and terminology that will hold throughout this paper.

Suppose $(Z,d_Z)$ is a path-metric space. If $A$ and $B$ are two subsets of $Z$, the \emph{distance between $A$ and $B$} is defined to be
$$
d_Z(A,B)=\inf\{d_Z(a,b)\ |\ a\in A, b\in B\}.
$$
A \emph{geodesic} in $Z$ is a path that minimizes distance globally. A geodesic segment between two points $x$ and $y$ is denoted $[xy]$. A \emph{geodesic triangle} $\Delta(xyz)$ is a union of three geodesic segments of the form $[xy]$, $[yz]$, and $[xz]$. A geodesic triangle $\Delta(xyz)$ in $Z$ is said to be \emph{$\delta$-slim} if for any point $w$ on $[xy]$ we have
$$
\min\{d_Z(w,[xz]), d_Z(w,[yz])\}\leq\delta.
$$

Suppose $Z$ has a cell-structure. The induced path-metric on the $k$--skeleton $Z^{(k)}$ of $Z$ is denoted $d_Z^k$. If $Y$ is any topological manifold, we will assume (by homotopy) that any map $f:Y\to Z$ induces via pullback a locally-finite cell structure on $Y$.

An \emph{edge-path} in $Z$ is a map $\sigma:[0,n]\to Z$ so that the restriction to each open interval $(i,i+1)$ is an isometry onto an open edge in $Z^{(1)}$. The \emph{length} of an edge-path $\sigma:[0,n]\to Z$, denoted $\ell(\sigma)$, is $n$. If $\sigma$ is an edge-path in $Z$ with $\sigma((i-1,i])\cap\sigma([i,i+1))=\sigma(i)$ for all $i$, then $\sigma$ is \emph{non-backtracking}. If $\sigma(0)=\sigma(n)$, then $\sigma$ is an \emph{edge-loop}. A \emph{minimal edge-loop} is an edge-loop minimizing length among other edge-loops in its free homotopy class in $Z$.

Throughout this paper, we let $N$ denote a closed irreducible (topological) 3--manifold with hyperbolic fundamental group $G$ containing a nontrivial closed surface subgroup $H\simeq\pi_1(S)$. Let $X$ denote the universal cover of $N$. Fix a 1--vertex triangulation of $N$ and let $d_N$ be the metric obtained by requiring all simplices of this triangulation to be isometric to standard regular Euclidean simplices. From this triangulation, one obtains a finite triangular presentation for $G$ with generators in 1--1 correspondence with the edges of the 1--skeleton $N^{(1)}$. This 1--skeleton lifts to $X^{(1)}$, which we identify with the Cayley graph for this presentation of $G$. The induced path-metric $d_X^1$ on $X^{(1)}$ is the standard Cayley-graph metric where all edges have length one.

The fact that the group $G$ is Gromov-hyperbolic is reflected in the fact that $(X^{(1)},d_X^1)$ has $\delta$-slim geodesic triangles, for some $\delta\geq 0$. The hyperbolicity of $G$ also implies that $S$ has genus $g\geq 2$.

Let $M$ denote the cover of $N$ with $\pi_1(M)=H\simeq\pi_1(S)$, and endow $M$ with the induced cell structure and metric. Equivalently, we have $M=X/H$, where $H$ acts in the usual way by left (say) translation. Note that $M$ is homotopy equivalent to $S$ (because closed hyperbolic surfaces are Eilenberg--Mac\ Lane spaces). We henceforth fix a homotopy equivalence $\phi:S\to M$, which provides a fixed isomorphism $\phi_*:\pi_1(S)\to H=\pi_1(M)$. This identification of $H$ as both a subgroup of $G$ and as the fundamental group of $S$ will often be tacitly assumed, especially with respect to notation. Suppose $\sigma$ is an essential closed curve on $S$. We let $\sigma^*$ denote a choice of minimal edge-loop in $M$ in the free homotopy class (in $M$) of $\phi(\sigma)$.

We say that $H$ is \emph{geometrically finite} in $G$ if there is some compact set $K\subset M$ so that, for all essential closed curves $\sigma$ on $S$, we have that $\sigma^*\subset K$. We say that $H$ is \emph{geometrically infinite} in $G$ if it is not geometrically finite.

Let $M_c$ be a compact core for $M$. An \emph{end} of $M$ is a non-precompact component of the complement $M-M_c$. A simple homological argument shows that $M$ has two ends. Note that if $K$ is any compact set containing $M_c$, there is a natural identification between the non-precompact components of $M-K$ with the ends of $M$. Such a component is called a \emph{neighborhood} of the corresponding end.

Thus $H$ is geometrically infinite in $G$ if and only if one may find a sequence $\{\sigma_i\}$ of essential curves on $S$ and a corresponding sequence of minimal edge-loops $\{\sigma_i^*\}$ in $M$ having the property that each neighborhood $U_b$ of some end $b$ intersects some $\sigma_i^*$ nontrivially. The sequence $\{\sigma_i^*\}$ is said to \emph{reach into} the end $b$. We say that $H$ is \emph{simply geometrically infinite} it is geometrically infinite, and if one may choose the $\sigma_i$ to be simple on $S$.

We now state our main theorem.
\begin{thm}\label{thm:main}
Let $N$ be a closed irreducible 3--manifold with fundamental group $G$ containing a closed surface subgroup $H$. The $H$ is geometrically infinite in $G$ if and only if $H$ is simply geometrically infinite in $G$.
\end{thm}

\begin{rem}
In the Kleinian groups context, one applies the terms ``geometrically finite'' and ``geometrically infinite'' to the \emph{ends} of the manifold $M$. It is more convenient in our context to apply these terms directly to the subgroup $H$. If one is thinking in terms of ends, however, our notion of simply geometrically infinite lies between those of a geometrically infinite end and a \emph{simply degenerate} end, for which one may choose the $\sigma_i$ not only to be simple on $S$ but also to have corresponding minimal edge-loops $\sigma_i^*$ that \emph{exit} the end (i.e., they leave compact sets while reaching deeper into the end). (In this case, we would say that $H$ is a \emph{simply degenerate subgroup} of $G$.) In this sense, the theorem above is a weaker version of the main result of~\cite{bonahon:ends}, which says essentially that geometrically infinite ends of hyperbolic 3--manifolds are simply degenerate. Here is a sketch of the argument that simply geometrically infinite implies simply degenerate in the hyperbolic 3--manifold case.

This is the pleated surface argument. Suppose $b$ is a simply geometrically infinite end of $M=\Hyp^3/H$ with $H\simeq\pi_1(S)$. Let $\{\sigma_i\}$ be a sequence of essential simple closed curves on $S$ so that the corresponding sequence of geodesics $\{\sigma_i^*\}$ in $M$ reaches into the end $b$ of $M$. Following an argument of Thurston, one can construct a sequence of pleated surfaces $\phi_i:S\to M$ so that $\phi_i(\sigma_i)=\sigma_i^*$. (Note that the existence of a map extending $\sigma_i\mapsto\sigma_i^*$ to all of $S$ depends crucially on the fact that $\sigma_i$ is simple on $S$.) Now each $\phi_i(S)$ is a union of totally geodesic triangles in $M$, with the number of triangles depending on the genus of $S$. Each of these triangles has area bounded above by $\pi$, so that the $\phi_i(S)$ have uniformly bounded area. Combining this fact with a uniform injectivity radius lower bound (in this set-up, this follows from the assumption that $M$ covers the compact manifold $N$), we deduce that each $\phi_i(S)$ has uniformly bounded diameter in $M$. In particular, the diameters of the $\sigma_i^*$ are uniformly bounded. The fact that they reach into $b$ therefore implies that they exit $b$. Thus $b$ is simply degenerate.

If one tries to trace through the argument just sketched assuming only that $\pi_1(N)$ is hyperbolic, the difficulty becomes obvious: geodesic triangles in arbitrary $\delta$-hyperbolic spaces need not have uniformly bounded areas. While there are a handful of obvious first steps around this difficulty, none have led so far to a proof. Nonetheless, it should be true that simply geometrically infinite surface subgroups of any torsion-free hyperbolic group are simply degenerate (the presence of an ambient 3--manifold should not be required for this).
\end{rem}

\subsection{Structure of the paper}\label{outline}

Here is an outline of the paper. In Section~\ref{intnum} we cover some basic results we need about $\delta$-hyperbolic spaces. In particular, we show that minimal edge-loops are coarsely unique, so that their coarse behavior depends only on the free homotopy class. We also prove the intersection number lemma. In Section~\ref{currents} we give some background on geodesic currents and laminations, and we fix the train-track notation that we will use in the remainder. In Section~\ref{shorten} we describe our version of Bonahon's splitting and straightening process, which is then used to prove Proposition~\ref{prop:ending}, which contains most of the heavy lifting in proving Theorem~\ref{thm:main}. This section is the most technical part of the paper, and also is where most of the novelty of this paper is to be found. Finally in Section~\ref{proof} we put the pieces together to prove Theorem~\ref{thm:main}.

As mentioned, both the basic outline of the proof and many of the detailed steps closely follow Bonahon's proof. Some of his techniques extend directly, and many are well-known, even in detail. As a result, some steps of our proof are presented as sketches, or simply with indications as to the changes required for our setting. To our knowledge, however, the details of Bonahon's intricate splitting and straightening argument appear only in the original (French) paper. Moreover, many of the details of this argument must be changed substantially in our setting. This section of the proof is therefore presented in full detail.

\subsection{Acknowledgements}\label{ack}

This paper began life as my doctoral thesis at the University of California, Santa Barbara. I am thus especially grateful to my thesis advisor Daryl Cooper for innumerable discussions and immeasurable help. I am also grateful to Sergio Fenley for pointing out an error in an earlier version of an attempt to show that simply geometrically infinite implies simply degenerate.

\section{Intersection Number}\label{intnum}

The following two results are well-known. We will prove the first, in part for the benefit of the reader unfamiliar with $\delta$-hyperbolic spaces, and in part because the techniques used in the proof will be invoked to justify a small step in the proof of Proposition~\ref{prop:ending}.

	\begin{lem}\label{lem:expgrowth}
	Let $\sigma^*$ be a minimal edge-loop in $M^{(1)}$, and suppose $\sigma$ is an edge-loop freely homotopic in $M$ to $\sigma^*$ with $d_M^1(\sigma,\sigma^*)=D$. Then there is a constant $c_1$, depending only on $\delta$, so that \[\ell(\sigma)\geq c_12^D\,\ell(\sigma^*).\]
	\end{lem}

	\begin{proof}
	The goal of the proof is to show that there is a coarse version of a nearest point retraction from $\sigma$ to $\sigma^*$ which is exponentially length increasing. The first step is to show that nearest point retraction is coarsely onto.

	Let $\tilde{\sigma}^*$ be a lift of $\sigma^*$ to $X$, and note that $\tilde{\sigma}^*$ is a geodesic in $X^{(1)}$. Now lift a homotopy between $\sigma^*$ and $\sigma$ to $X$ to obtain a corresponding lift $\tilde{\sigma}$ of $\sigma$.

Let $T\subset X^{(0)}$ be those vertices $t\in\tilde{\sigma}^*$ so that there exists a vertex $s\in\tilde{\sigma}$ with $d_X^1(s,\tilde{\sigma}^*)=d_X^1(s,t)$. We will show that the components of $\tilde{\sigma}^*-T$ are no longer than $7\delta$.

	Choose two adjacent vertices $s$ and $s'$ on $\tilde{\sigma}$, and suppose $t$ and $t'$ are vertices on $\tilde{\sigma}^*$ nearest to $s$ and $s'$, respectively (see Figure~\ref{coarseonto}). Suppose $d_X^1(t,t')\geq 7\delta$. Then we may choose a vertex $z_0$ on $[tt']$ with $d_X^1(z_0,t)>2\delta$ and $d_X^1(z_0,t')>4\delta$. Because the geodesic triangle $\Delta(tst')$ is $\delta$-slim, there is a vertex $z_1\in[ts]\cup [st']$ with $d_X^1(z_0,z_1)\leq\delta$. If $z_1$ were on $[ts]$, however, then joining $z_0$ to $z_1$ and $z_1$ to $t$ would produce a path with length no more than $2\delta$ joining $z_0$ to $t$, contradicting that $\tilde{\sigma}^*$ is a minimal edge-path.

	\begin{figure}[ht]
	\begin{center}
		\psfrag{s}{$s$}
		\psfrag{r}{$s'$}
		\psfrag{t}{$t$}
		\psfrag{q}{$t'$}
		\psfrag{A}{$\tilde{\sigma}^*$}
		\psfrag{B}{$\tilde{\sigma}$}
		\psfrag{z}{$z_0$}
		\psfrag{w}{$z_1$}
		\includegraphics{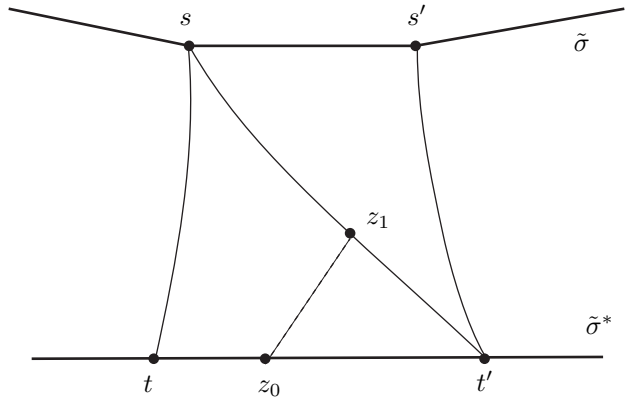}
	\end{center}
	\caption{Nearest point retraction is coarsely onto}
	\label{coarseonto}
	\end{figure}

	Thus $z_1$ is contained on $[st']$. Now the triangle $\Delta(t'ss')$ is $\delta$-slim, so we find $z_2\in[ss']\cup[s't']$ with $d_X^1(z_1,z_2)\leq\delta$. As above, we deduce from the minimality of $\tilde{\sigma}^*$ and the choice of $z_0$ that $z_2$ lies on $[ss']$. But this implies that either $z_2=s$ or $z_2=s'$, each of which contradicts the assumption that $d_X^1(\sigma,\sigma^*)\geq D>2\delta$. Thus the components of $\tilde{\sigma}^*-T$ are no longer than $7\delta$.

	Now fix a vertex $\tilde{p}_0\in T\subset\tilde{\sigma}^*$, and let $\tilde{p}_0'$ be a vertex on $\tilde{\sigma}^*$ with $d_X^1(\tilde{p}_0,\tilde{p}_0')=12\delta$. Then there is some vertex $\tilde{p}_1\in T$ with $d_X^1(\tilde{p}_0',\tilde{p}_1)\leq 4\delta$, and so in particular, we have $d_X^1(\tilde{p}_0,\tilde{p}_1)\geq 8\delta$. Let $\tilde{q}_0$ be a vertex on $\tilde{\sigma}$ with $d_X^1(\tilde{q}_0,\tilde{p}_0)=d_X^1(\tilde{q}_0,\tilde{\sigma}^*)$, and similarly choose $\tilde{q}_1$ corresponding to $\tilde{p}_1$.

	Consider the ``quadrilateral'' in $X^{(1)}$ formed by that portion of $\tilde{\sigma}^*$ between $\tilde{p}_0$ and $\tilde{p}_1$, that portion of $\tilde{\sigma}$ between $\tilde{q}_0$ and $\tilde{q}_1$, and minimal edge-paths $[\tilde{p}_0\tilde{q}_0]$ and $[\tilde{p}_1\tilde{q}_1]$ (see Figure~\ref{midpoints}). Let $L$ be the length of the portion of $\tilde{\sigma}$ between $\tilde{q}_0$ and $\tilde{q}_1$.

	\begin{figure}[ht]
	\begin{center}
		\psfrag{s}{$\tilde{\sigma}$}
		\psfrag{S}{$\tilde{\sigma}^*$}
		\psfrag{p}{$\tilde{p}_i$}
		\psfrag{P}{$\tilde{p}_{i+1}$}
		\psfrag{q}{$\tilde{q}_i$}
		\psfrag{Q}{$\tilde{q}_{i+1}$}
		\psfrag{m}{$m$}
		\psfrag{l}{$m_0$}
		\psfrag{n}{$m_1$}
		\psfrag{0}{$z_0$}
		\psfrag{1}{$z_1$}
		\psfrag{2}{$z_2$}
		\psfrag{t}{$z_t$}
		\psfrag{i}{$k_{10}$}
		\psfrag{j}{$k_{11}$}
		\psfrag{k}{$k_1$}
		\psfrag{K}{$k$}
		\includegraphics{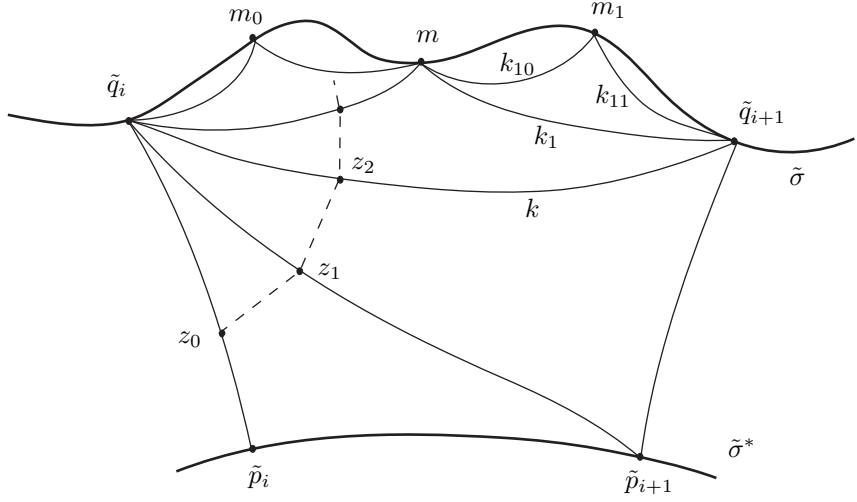}
	\end{center}
	\caption{Exponential growth of curves}
	\label{midpoints}
	\end{figure}

	Because $D>2\delta$, there is a vertex $z_0\in[\tilde{p}_0\tilde{q}_0]$ with $\delta<d_X^1(z_0,[\tilde{p}_0\tilde{p}_1])<2\delta$. Using arguments similar to those above, there is a vertex $z_1\in[\tilde{p}_1\tilde{q}_0]$ with $d(z_0,z_1)\leq\delta$. Similarly, there is a vertex $z_2$ on $[\tilde{q}_0\tilde{q}_1]\cup[\tilde{q}_1\tilde{p}_1]$ with $d_X^1(z_1,z_2)\leq\delta$. If $z_2\in[\tilde{q}_1\tilde{p}_1]$, however, then \[d_X^1(\tilde{p}_0,\tilde{p}_1)\leq d_X^1(\tilde{p}_0,z_0)+d_X^1(z_0,z_1)+d_X^1(z_1,z_2)+d_X^1(z_2,\tilde{p}_1)\] \[<4 \delta+d_X^1(z_2,[\tilde{p}_0\tilde{p}_1])\leq 8\delta,\] contradicting the assumption that $d_X^1(\tilde{p}_0,\tilde{p}_1)\geq 8\delta$. It follows that $z_2\in[\tilde{q}_0\tilde{q}_1]$.

	In order now to show that the length of that part of $\tilde{\sigma}$ between $\tilde{q}_0$ and $\tilde{q}_1$ depends exponentially on $D$, we mimic the procedure used in the proof of~\cite[Theorem $2.19$]{short}. By a \emph{midpoint} of a subarc of $\tilde{\sigma}$ we mean the vertex on $\tilde{\sigma}$ which divides the subarc into two pieces whose lengths differ by no more than one. We let $k$ be a minimal edge-path between $\tilde{q}_0$ and $\tilde{q}_1$. Let $m$ be the midpoint of the arc of $\tilde{\sigma}$ between $\tilde{q}_i$ and $\tilde{q}_1$, with $k_0$ a minimal edge-path from $\tilde{q}_0$ to $m$ and $k_1$ a minimal edge-path between $m$ and $\tilde{q}_1$.

	Now let $b$ be a binary sequence and suppose $k_b$ to have been chosen, joining $k_b(0)$ to $k_b(1)$. Let $m_b$ be the midpoint of the portion of $\tilde{\sigma}$ between the endpoints of $k_b$. Then we let $k_{b0}$ denote a minimal edge-path between $k_b(0)$ and $m_b$, while $k_{b1}$ denotes a minimal edge-path between $m_b$ and $k_{b1}$. We continue this construction until all $k_b$ have length one, and thus are edges in $\tilde{\sigma}$.

	Note that there is some vertex $z_3$ on either $k_0$ or $k_1$ with $d_X^1(z_2,z_3)\leq\delta$. More generally, if $z_j$ lies on $k_b$, there is some vertex $z_{j+1}$ on either $k_{b0}$ or $k_{b1}$ with $d_X^1(z_j,z_{j+1})\leq\delta$. The $z_j$ then provide a path from $z_0\in\tilde{\sigma}^*$ to $z_t\in\tilde{\sigma}$, where $t\leq\log_2(L)+2$. Thus the length of the path is no more than $\delta(\log_2(L)+2)$. On the other hand, this path must have length at least $D$, so we deduce that \[D\leq\delta(\log_2(L)+2),\] and so \[L\geq\frac{2^{D/\delta}}{4}.\]

	Now let $n$ be the greatest integer in $\ell(\sigma^*)/12\delta$. If $n\geq 1$, then by repeating this analysis on $n$ segments of $\sigma^*$ of length $12\delta$, we find that $\ell(\sigma)\geq n\frac{2^{D/\delta}}{4}$. But $\ell(\sigma^*)/12\delta\leq n+1$, and then because $n\geq 1$ we have $(n+1)/2\leq n$, so $\ell(\sigma^*)/24\delta<n$. Thus we obtain \[\ell(\sigma)\geq\frac{2^{D/\delta}}{96\delta}\,\ell(\sigma^*).\]

	If, on the other hand, we have that $n=0$, then we set $m$ to be the greatest integer in $12\delta/\ell(\sigma^*)$. Then the segment of $\tilde{\sigma}^*$ between $\tilde{p}_0$ and $\tilde{p}_1$ contains an $m$-fold lift of $\sigma^*$. Thus each lift contributes to at least a portion of $1/m$ to the expansion of $\ell(\sigma^*)$ to $L\geq\frac{2^{D/\delta}}{4}$. In other words, we have $\ell(\sigma)\geq\frac{1}{m}\frac{2^{D/\delta}}{4}$. But then $\ell(\sigma^*)/12\delta\leq\frac{1}{m}$, and so we have \[\ell(\sigma)\geq\frac{2^{D/\delta}}{48\delta}\ell(\sigma^*).\]
	\end{proof}

The second result we require is the following (cf.~\cite[Lemma $2.1$]{bonahon:ends}).

	\begin{lem}\label{lem:nbhd}
	Let $\sigma^*\subset M^{(1)}$ be a minimal edge-loop and suppose $\sigma$ is an edge-loop freely homotopic to $\sigma^*$ in $M$ formed of $n$ geodesic segments in $M^{(1)}$. Then there is a constant $c_2$, depending only on $\delta$, so that $\sigma^*$ is contained in the $nc_2$--neighborhood of $\sigma$.
	\end{lem}

Note that this lemma implies (by choosing $n=1$) that the choice of minimal edge-loop $\sigma^*$ corresponding to $\sigma$ described above is coarsely unique, in the sense that any two such choices must have Hausdorff distance bounded above by $c_2$.

We also require the following annular version of the linear isoperimetric inequality (cf.~\cite[\S $7.4.$B]{gromov:hyperbolic}).
	\begin{lem}\label{lem:isoannulus}
	There is a constant $K$, depending only on $\delta\geq 1$, with the following property: if $\sigma$ is an essential edge-loop in $M$, $\sigma^*$ is a corresponding minimal edge-loop, and $A\co S^1\times[0,1]\to M^{(2)}$ is a minimal-area annulus with boundary $\sigma\cup\sigma^*$, then \[\text{area}(A)\leq K(\ell(\sigma)+\ell(\sigma^*)).\]
	\end{lem}

The first application of the preceeding results is the proof of the intersection number lemma, Lemma~\ref{lem:intnum} below. Our statement is the same as that of~\cite[Prop $3.4$]{bonahon:ends}, and the proof is actually simplified somewhat by the combinatorial setting.

	\begin{lem}\label{lem:intnum}
	Let $\alpha_1^*$ and $\alpha_2^*$ be two minimal edge-loops of $M$ homotopic to two curves $\alpha_1$ and $\alpha_2$ on $S$ by two homotopies not otherwise hitting $S$ and arriving on the same side of $S$. Then there exists a constant $C$, depending only on $\delta$, so that if $d_M^1(\alpha_i^*,S)\geq D$ for $i=1,2$, then
	\[\mathrm{int}(\alpha_1,\alpha_2)\leq C\ell(\alpha_1)\ell(\alpha_2)2^{-D},\] where $\mathrm{int}(\ ,\ )$ and $\ell(\ )$ denote respectively intersection number in $S$ and length in $M$.
	\end{lem}

	\begin{proof}
	Let $A_2\co S^1\times[0,1]\to M^{(2)}$ be an annular homotopy between $\alpha_2$ and $\alpha_2^*$ which minimizes area in its homotopy class in $M$ (relative to its boundary). We abuse notation in the usual way by letting $A_2$ denote the image of the annulus in $M$, just as we let $\alpha_1^*$ denote the image of the edge-path in $M$.

	Note that by pushing the interiors of the 1--cells of $\alpha_1^*$ into the interiors of the 3--cells of $M$, while fixing the vertices, we may ensure that all intersections of $\alpha_1^*$ with $A_2$ lie on vertices in $M$. Thus the geometric intersection number of $A_2$ with $\alpha_1^*$ is bounded above by the number of vertices in the abstract product $S^1\times(S^1\times[0,1])$ of the domains (with the induced combinatorial structures). But the number of vertices in $A_2$ is bounded above by \[\frac{\text{area}(A_2)}{2}+\ell(\alpha_2)+\ell(\alpha_2^*)\leq\frac{\text{area}(A_2)}{2}+2\ell(\alpha_2).\] Also we have from Lemma~\ref{lem:isoannulus} that \[\text{area}(A_2)\leq K\big(\ell(\alpha_2^*)+\ell(\alpha_2)\big)\leq 2K\ell(\alpha_2).\] Thus we have \[\text{int}(\alpha_1^*,A_2)\leq\ell(\alpha_1^*)\text{area}(A_2)\leq(K+2)\ell(\alpha_1^*)\ell(\alpha_2).\]

On the other hand, we have from Lemma~\ref{lem:expgrowth} that \[\ell(\alpha_1^*)\leq c_1^{-1}2^{-D}\ell(\alpha_1).\] The preceding two inequalities together imply that \[\text{int}(\alpha_1^*,A_2)\leq (K+2)c_1^{-1}\ell(\alpha_1)\ell(\alpha_2)\,2^{-D}.\] By symmetry, $\text{int}(\alpha_2^*,A_1)$ is bounded in the same way.

We now use the fact (cf.~\cite{bonahon:ends}, Lemma $3.2$) that \[\text{int}(\alpha_1,\alpha_2)\leq\text{int}(\alpha_1^*,A_2)+\text{int}(\alpha_2^*,A_1)\] to deduce that \[\text{int}(\alpha_1,\alpha_2)\leq 2(K+2)c_1^{-1}\ell(\alpha_1)\ell(\alpha_2)\,2^{-D}.\]
	\end{proof}

\section{Geodesic Currents \& Train Tracks}\label{currents}

In this section we discuss briefly the tools necessary to apply Bonahon's lamination-straightening argument. Other than some minor adjustments to account for the change is setting, this is primarily a review of the material in~\cite[\S\S $3.1$--$5.1$]{bonahon:ends}.

\subsection{Geodesic Currents}\label{}

A concise account of the necessary facts about geodesic currents may be found in~\cite[\S IV]{bonahon:ends}, on which the following discussion is based. More details may be found in~\cite{bonahon:currents}.

Given a closed Riemannian surface $S$ we consider the projective tangent bundle $PT(S)$, which is the quotient of the unit tangent bundle \[UT(S)=\{(p,v)\ \mid\ p\in S,\ v\in T_p(S),\ \norm{v}=1\}\] obtained by identifying $(p,v)$ and $(p,-v)$. Thus a point in $PT(S)$ is a pair $(p,[v])$, where $[v]=\{v,-v\}$. A point $(p,v)\in UT(S)$ uniquely determines a unit speed geodesic $\gamma(t)$ on $S$ with $\gamma(0)=p$ and $\gamma'(0)=v$. We define a flow $\Phi$ on $UT(S)$, called the \emph{geodesic flow}, by setting $\Phi((p,v),t)=(\gamma(t),\gamma'(t))$. The quotient of the flow lines under the involution sending $v\in T_p(S)$ to $-v$ forms a 1--dimensional foliation $\mathscr{F}$ of $PT(S)$.

A \emph{geodesic current} on a surface $S$ is a positive transverse measure $\mu$ invariant under the geodesic foliation $\mathscr{F}$. In other words, given a codimension--1 submanifold $V$ of $PT(S)$ which is transverse to the foliation, we require that $\mu(V)=\mu(\Phi(V,t))$ for all $t$. For the most basic example, suppose $\alpha$ is a closed (not necessarily simple) geodesic on $S$. Then $\alpha$ lifts to a leaf $\tilde{\alpha}$ of $PT(S)$. We may thus define a geodesic current, also denoted by $\alpha$, by declaring that $\alpha(V)$ equals the number of times $V$ intersects the lift $\tilde{\alpha}$ of $\alpha$ to $PT(S)$.

\begin{figure}[ht]
	\begin{center}
		\includegraphics{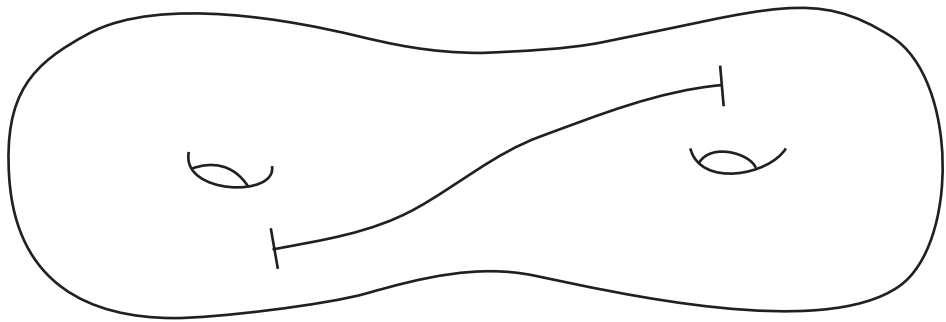}
	\end{center}
	\caption{An elongated $H$ defining a flow box.}
	\label{box}
	\end{figure}

We let $\mathscr{C}(S)$ denote the space of geodesic currents on $S$ with topology defined as follows. Consider an elongated letter $H$ on $S$, with the horizontal bar a geodesic arc transverse to the vertical bars, and with the vertical bars short enough so that all geodesic arcs joining the two vertical bars and homotopic to a path in the $H$ intersect the vertical bars transversely (see Figure~\ref{box}). A \emph{flow box} $B\subset PT(S)$ for $\mathscr{F}$ is formed by lifting all geodesic arcs on $S$ joining the vertical bars of the $H$ and homotopic to a path in $H$. If we let $Q$ denote the abstract product of the two vertical arcs of $H$, then there is a diffeomorphism $B\cong Q\times[0,1]$ so that the leaves of $B\cap\mathscr{F}$ correspond to $(\text{point})\times[0,1]$. We let $\partial_{\mathscr{F}}B$ denote that portion of $B$ corresponding to $\partial Q\times[0,1]$.

Given a geodesic current $\alpha\in\mathscr{C}(S)$ and a flow box $B$, we define the measure $\alpha(B)\in\R^+$ to be the measure under $\alpha$ of $Q\times(\text{point})$ under the identification $B\cong Q\times[0,1]$. Thus if $\alpha$ is a geodesic current corresponding to a closed geodesic on $S$ as above, then $\alpha(B)$ is exactly the number of subarcs of $\alpha$ which lift to leaves of $B\cap\mathscr{F}$.

A neighborhood basis of $\alpha\in\mathscr{C}(S)$ is formed of open sets in $\mathscr{C}(S)$ of the form \[u(\alpha;B_1,\dots,B_n;\epsilon)=\{\beta\in\mathscr{C}(S)\ \mid\ \forall i, |\alpha(B_i)-\beta(B_i)|<\epsilon\},\] where $\epsilon>0$ and the $B_i$ are a collection of flow boxes so that $\alpha(\partial_{\mathscr{F}}(B))=0$.

Geodesic currents are a generalization of measured geodesic laminations (see \cite{casson-bleiler,travaux}). One may think of geodesic currents as the non-simple analogues of measured geodesic laminations. In particular, Bonahon shows that the intersection number function defined on the set of closed geodesics on $S$ extends to a continuous bilinear function $\text{int}\co\mathscr{C}(S)\times\mathscr{C}(S)\to\R^+$, and that geodesic measured laminations are exactly those geodesic currents with zero self-intersection number.

The length function defined on closed weighted geodesics on $S$ also extends to a linear continuous non-negative function on the space of geodesic currents. What is of interest to us is the projective class of a current (in $\mathcal{PC}(S)$), so we will work primarily with unit-length currents. Given a closed geodesic $\gamma$ on $S$, we let $\overline{\gamma}$ denote the unit length current represented geometrically by $\gamma$, and we write $\overline{\gamma}=\gamma/\ell(\gamma)$.

It is possible to put a metric on $\mathscr{PC}(S)$ making it into a complete space containing collections of closed (projective) weighted geodesics as a dense subset. In particular, after subsequencing, any infinite sequence $\{\gamma_i\}$ of closed geodesics on $S$ gives rise to a convergent sequence $\{\overline{\gamma}_i\}$ of unit-length geodesic currents.

\subsection{Train tracks}\label{tracks}

\begin{figure}[ht]
	\begin{center}
		\psfrag{c}{corner}
		\psfrag{t}{tie}
		\psfrag{s}{switch}
		\includegraphics{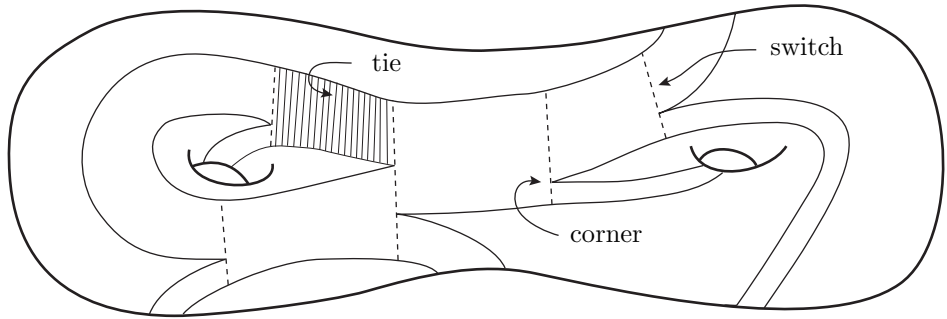}
	\end{center}
	\caption{A train track on a surface}
	\label{traintrack}
	\end{figure}

Following~\cite{bonahon:ends}, we define a \emph{train track} $\tau$ on a surface $S$ to be a closed subset of $S$ which can be written as a finite union of rectangles $R_i$ with disjoint interiors, each of which carries a vertical foliation by arcs, called \emph{ties}. These rectangles are glued together along their vertical (``short'') boundaries $\partial_vR_i$ so that the boundary of $\tau$ is the union of the horizontal (``long'') boundaries $\partial_hR_i$ (see Figure~\ref{traintrack}). We also insist that the boundary of $\tau$ be smooth except at a finite number of points, called \emph{corners}, where $\partial_hR_i$ meets the interior of $\partial_vR_j$ for some $i\neq j$. Finally, we will assume that no connected component of $-\tau$ has closure a disk with two or fewer corners, or an annulus with no corners.

The vertical components of the boundaries of the rectangles are called the \emph{switches} of $\tau$, and the \emph{branches} are the closures of the connected components of the complement of the switches.

A \emph{train path} is a path in the interior of a train track $\tau$ which is transverse to the ties. We say that $\alpha$ is \emph{carried by} $\tau$ if $\alpha$ is isotopic to a train path. If $\alpha$ is isotopic to a train path intersecting all ties of $\tau$, then we say that $\alpha$ is \emph{fully carried} by $\tau$. Every lamination $\alpha$ on $S$ is fully carried by some train track, as, for $\epsilon>0$ small enough, the $\epsilon$--neighborhood of $\alpha$ on $S$ is a traintrack.

For a measured lamination $\alpha$ fully carried by $\tau$, we obtain a \emph{weighted train track} where the \emph{weight} of a branch $b$ is exactly the measure $\alpha(b)$ under $\alpha$ of one of its ties. Note that $\tau$ may have annular components corresponding to compact leaves of $\alpha$.

Given a train track $\tau$ fully carrying a lamination $\alpha$, Bonahon shows that we may cover $PT(S)$ by a finite collection of flow boxes $B_1,\dots,B_q$ with disjoint interiors so that for $1\leq i\leq p$, the arcs of $B_i\cap\mathscr{F}$ project into $S$ as train paths in $\tau$ properly embedded in the branches, while the lift of $\alpha$ to $PT(S)$ is disjoint from the boxes $B_j$ for $p<j\leq q$. To do this, we first cover the support of $\alpha$ by boxes $B_1,\dots,B_p$ with the required property. These are defined by elongated $H$'s whose horizontal bars are contained in the support of $\alpha$ and whose vertical bars are (non-overlapping) portions of switches of $\tau$. We extend this collection to a cover of $PT(S)$ by adding boxes $B_{p+1}'\dots,B_r'$ so that all arcs of $B_j'\cap\mathscr{F}$ project to arcs disjoint from $\alpha$. We further suppose that the vertical bars of the $H$'s defining the $B_j'$ are disjoint from the $H$'s defining the $B_i$ with $i\leq p$, and are pairwise disjoint. We finally observe that each of $B_j'\cap B_k'$, $\overline{B_j'-B_k'}$ and $\overline{B_j'-B_i}$ is the union of a finite number of flow boxes with disjoint interiors.

Thus if $b$ is a branch of a train track $\tau$ carrying a lamination $\alpha$, the measure $\alpha(b)$ is equal to the sum of the measures under $\alpha$ of those boxes of the $B_1,\dots,B_p$ the arcs of which project to train paths in $b$.

We will need the following standard fact about train tracks, which follows immediately from the lack of disks, monogons, digons, or annuli in the complement of $\tau$.
	\begin{lem}\label{lem:track}
	In a train track $\tau$ on a surface $S$, two curves transverse to the ties and homotopic in $S$ (with endpoints fixed) are homotopic in $\tau$ by a homotopy preserving the ties.
	\end{lem}

\begin{figure}[ht]
	\begin{center}
		\psfrag{A}{$\downarrow$}
		\psfrag{t}{$\tau$}
		\psfrag{T}{$\tau'$}
		\includegraphics{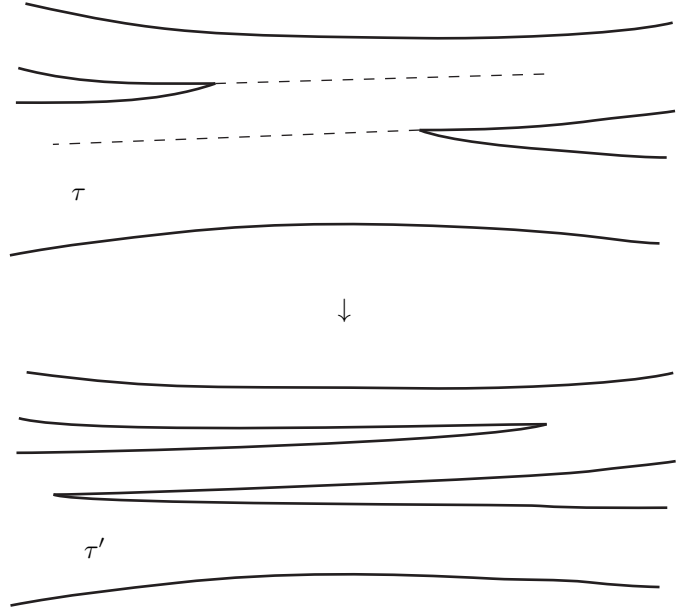}
	\end{center}
	\caption{Splitting a train track}
	\label{splittrack}
	\end{figure}

A \emph{splitting} of $\tau$ on $S$ is a train track $\tau'\subset\tau$ whose ties are contained in ties of $\tau$. Note that any splitting of $\tau$ may be further split so as to have the same number of branches as $\tau$. Also note that all laminations carried by any splitting of $\tau$ are carried by $\tau$. Given a train track $\tau$ carrying a lamination $\alpha$, we may obtain a splitting $\tau'$ of $\tau$ which also carries $\alpha$ by cutting along a finite collection of pairwise disjoint arcs in $\tau$ which are disjoint from $\alpha$ and transverse to the ties, and each of which has one endpoint on a corner of $\tau$ (see Figure~\ref{splittrack}).
Note that there is a bound on the possible number of annular components of any train track $\tau$ which fully carries a lamination $\alpha$ on  $S$. In particular, any such train track has a splitting which already has as many annular components as any further splitting.

\section{Shortening Laminations}\label{shorten}

In perhaps the most technical section of~\cite{bonahon:ends}, a train track $\tau$ fully carrying a geodesic measured lamination $\alpha$ is repeatedly split, and these splittings are mapped into the 3--manifold, where the length of its image can be compared to its length on the surface. By carefully choosing the splittings of $\tau$ and the corresponding maps, one deduces (\cite[Thm $5.1$]{bonahon:ends}) that either the length of $\phi(\alpha)$ is tending to zero, or $\phi(\alpha)$ is becoming very straight. More precisely, if for any map $\phi:S\to M$ we have that the length of $\phi'(\alpha)$ over all $\phi'$ homotopic to $\phi$ is bounded away from zero, then by proceeding far enough along in the straightening process, one may ensure that closed curves sufficiently close (as currents) to $\alpha$ have images which run along their geodesic representatives at an arbitrarily small distance for an arbitrarily large percentage of their length.

There is a corresponding, although slightly weaker, statement to be made in this context, the proof of which is similar to that of~\cite[Thm $5.1$]{bonahon:ends}. Recall that, given an essential closed curve $\gamma$ on $S$, $\gamma^*$ denotes a (coarsely unique) choice of minimal edge-loop freely homotopic in $M$ to $\phi(\gamma)$, where $\phi:S\to M$ is a fixed homotopy equivalence.
	\begin{prop}\label{prop:ending}
	Let $\alpha$ be a unit length geodesic measured lamination on $S$. Fix a number $t<1$. There is a constant $\epsilon$, depending only on $t$ and the hyperbolicity constant $\delta$, so that either:
	\begin{enumerate}
	\item there exists a map $\phi'\co S\to M$ homotopic to $\phi$ so that for any closed curve $\gamma$ on $S$, the ratio, to the total length of $\phi'(\gamma)$, of the length of that portion of $\phi'(\gamma)$ running along $\gamma^*$ at a distance less than $\epsilon$, is greater than $t$ when $\overline{\gamma}$ is sufficiently close to $\alpha$ in $\mathscr{PC}(S)$, or
	\item we may homotope $\phi$ to make $\ell(\phi(\alpha))$ arbitrarily small.
	\end{enumerate}
	\end{prop}

In the following sections, we first discuss the types of maps we will be considering, and explain exactly what is meant by the length $\ell(\phi(\alpha))$. We then define the splitting procedure and identify its salient properties. Finally we prove the proposition by repeated application of the procedure.

This conclusion of this theorem is essentially the same as that used by Bonahon except that $\epsilon$ depends on $t$. Given that minimal edge-paths are not unique (and, in particular, do not always have Hausdorff distance zero), one does not expect to be able to vary $\epsilon$ freely.

\begin{rem}
As mentioned in the introduction, the proof of this proposition makes no use of the fact that $G$ is assumed to be a 3--manifold group. In particular, given a torsion-free hyperbolic group $G$, we may take $N$ to be a finite $K(G,1)$, and let $M$ be the cover of $N$ corresponding to $H=\pi_1(S)$. The conclusion of the proposition will hold in this case.
\end{rem}

\subsection{Adapted Maps}\label{maps}

We say a map $\phi\co S\to M^{(2)}$ is \emph{adapted} to a train track $\tau$ if it sends each branch of $\tau$ to a non-backtracking edge-path in $M$, so that the preimage of a point $x\in\phi(\tau)$ is a union of ties, and so that $\phi$ is locally injective on any path in $\tau$ transverse to the ties. We say that $\phi$ is \emph{tightly adapted} to $\tau$ if it is adapted to $\tau$ so that each branch is sent to a minimal edge-path in $M$. In particular, if $\phi$ is tightly adapted to $\tau$, then $\phi$ maps any compact leaves of $\alpha$ carried by annular components of $\tau$ to minimal edge-loops in $M$.

Suppose $\alpha$ is a unit length geodesic measured lamination on $S$ carried by the train track $\tau$, and suppose $\phi\co S\to M$ is a map adapted to $\tau$. Pulling the path metric on $M^{(1)}$ back along $\phi$, one obtains measures for the lengths of the branches of $\tau$. With the usual abuse of notation, we let $\ell(\phi(b_i))$ denote the length of a branch $b_i$ measured in this way, and define the \emph{length of $\phi(\alpha)$ in $M$} to be \[\ell(\phi(\alpha))=\sum_i\alpha(b_i)\ell(\phi(b_i)),\] where the sum is over all branches $b_i$ of $\tau$ and $\alpha(b_i)$ is the measure under $\alpha$ of a tie of $b_i$.

We define a \emph{vertex tie} for the pair $(\phi,\tau)$ to be any tie $t$ of $\tau$ with $\phi(t)$ a vertex in $M$. Now let $\phi\co S\to M$ be adapted to a train track $\tau$ and suppose $\tau'$ is a splitting of $\tau$ on $S$, chosen so that each switch of $\tau'$ is contained in a vertex tie of $\tau$. Then we may construct a new map $\phi'$ tightly adapted to $\tau'$ as follows: For each switch $s$ of $\tau'$ we define $\phi'(s)=\phi(s)$. Now for each branch $b$ in $\tau'$ we require that $\phi'$ map $b$ to a minimal edge-path in $M^{(1)}$ homotopic with fixed boundary to $\phi(b)$. We then extend arbitrarily, mapping the remainder of $S$ into $M^{(2)}$. The fact that $\phi$ and $\phi'$ are homotopic follows from the irreducibility of $M$. With $\phi'$ obtained in this way, we clearly then have that \[\ell(\phi'(\alpha))\leq\ell(\phi(\alpha)).\]

\subsection{The splitting procedure}\label{splitting}

Let $\phi\co S\to M$ be a map adapted to a train track $\tau$, and fix $\epsilon>0$. We define an \emph{$\epsilon$-shortcut} for the pair $(\tau,\phi)$ to be an arc $k\subset\tau$ with the following properties:
	\begin{enumerate}
	\item the endpoints of $k$ are contained in vertex ties,
	\item $k$ is transverse to the ties
	\item $\phi(k)$ is a non-minimal edge-path in $M$ homotopic with endpoints fixed to an edge-path of length no more than $\epsilon$.
	\end{enumerate}
We consider two shortcuts equivalent if they are homotopic by a tie-preserving homotopy. For simplicity, we will suppress reference to equivalence classes and speak only of shortcuts.

For any subset $X$ of a train track $\tau$, we let $R(X)$ denote the union of those ties of $\tau$ which meet $X$. When $\tau$ is adorned with sub- or superscripts, we similarly adorn $R$. Thus, for example, $R'(X)$ is the union of ties of $\tau'$ meeting $X$.

	\begin{lem}(cf.~\cite[Lemma $5.7$]{bonahon:ends})\label{lem:shortcuts}
	Suppose $\tau$ is a train track carrying a lamination $\alpha$. Then there is a splitting $\tau'$ of $\tau$, still carrying $\alpha$, and a finite collection $\{k_{\mu},\dots,k_{\nu}\}$ of $\epsilon$-shortcuts for $\phi(\tau')$ with the following properties:
	\begin{enumerate}
	\item $R'(k_i)$ and $R'(k_j)$ have disjoint interiors, for $i\neq j$;
	\item any other $\epsilon$-shortcut for $\phi(\tau')$ has at least one endpoint in the interior of some $R'(k_i)$;
	\item each $k_i$ crosses at most one switch of $\tau'$.
	\end{enumerate}
	\end{lem}

Before beginning the proof, we note that although the train track $\tau$ is being split to obtain $\tau'$, the map $\phi$ is unchanged. It follows from the definitions that if $\phi$ was adapted to $\tau$ then $\phi$ is adapted to any splitting $\tau'$ of $\tau$. Property (3) in the lemma allows for the construction of $\phi'$, also adapted to $\tau'$, by straightening the shortcuts $k_{\mu},\dots,k_{\nu}$. Properties (1) and (2) allow for an estimation in the decrease in length in passing from $\phi(\tau)$ to $\phi'(\tau')$.

	\begin{proof}
	We note that the proof that follows is essentially the same as that of~\cite[Lemma $5.7$]{bonahon:ends}. To begin, we may assume that $\tau$ contains as many annular components as any further splitting, and that $\phi$ is tightly adapted to $\tau$. In particular, any annular components of $\tau$ are mapped to minimal edge-loops in $M$, and thus contain no shortcuts.

	Note that there are only finitely many homotopy classes of paths shorter than $\epsilon$ in $M$ that join the images under $\phi$ of any two branches of $\tau$. Because $\phi_*:\pi_1(S)\to\pi_1(M)$ is injective, there are only finitely many $\epsilon$-shortcuts in $S$ joining any two particular branches. In particular, there is some constant $A$ (depending on $\phi$) with the property that any $\epsilon$-shortcut $k$ for $\phi(\tau)$ has $\ell(\phi(k))\leq A$.

	By the lemma below, we may split $\tau$ to $\tau_0$ so that the non-annular branches of $\phi(\tau_0)$ all have length in $M$ greater than $2A$. Thus an $\epsilon$-shortcut for $\phi(\tau_0)$ hits at most one switch of $\tau_0$ (but note that because we have not altered the map $\phi$, some $\epsilon$-shortcuts may not cross any switches).

	We first find a finite family $\{k_1,\dots,k_p\}$ of $\epsilon$-shortcuts not crossing any switches of $\tau_0$ chosen to be maximal subject to the condition that the rectangles $R_0(k_i)$ have disjoint interiors (see Figure~\ref{shortcuts}). Note that all other $\epsilon$-shortcuts not crossing any switches of $\tau_0$ necessarily have an endpoint in the interior of some $R_0(k_i)$, as otherwise there is a collection of shortcuts including this one which contradicts the maximality condition.

	\begin{figure}[ht]
	\begin{center}
		\psfrag{t}{$\tau_0$}
		\psfrag{i}{$k_i$}
		\psfrag{j}{$k_j$}
		\psfrag{U}{$U_0$}
		\psfrag{o}{other shortcuts}
		\includegraphics{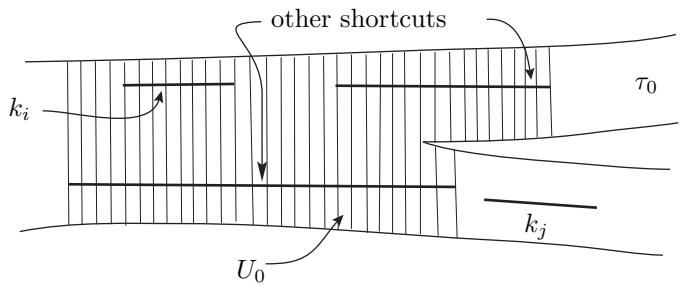}
	\end{center}
	\caption{Finding a maximal family of shortcuts}
	\label{shortcuts}
	\end{figure}

	Let $U_0$ be the union of the $R_0(k)$ where $k$ is an $\epsilon$-shortcut with endpoints not contained in any of the $R_0(k_i)$ (in particular, any such $k$ must cross a switch of $\tau_0$). Note that because the branches of $\tau_0$ have length greater than $2A$, we may join each corner of $\tau_0$ contained in $U_0$ to a tie in $\partial U_0$ by a path contained in $U_0$, transverse to the ties, not crossing any switches of $\tau_0$, and disjoint from $\alpha$. We choose these paths to be pairwise disjoint, and then split along them (see Figure~\ref{splittingshortcuts}).

	\begin{figure}[ht]
	\begin{center}
		\psfrag{t}{$\tau'$}
		\psfrag{j}{$k_j,\ \mu\leq j\leq p$}
		\psfrag{s}{$k_s,\ s>p$}
		\psfrag{i}{$k_i,\ i<\mu$}
		\psfrag{U}{$U'$}
		\includegraphics{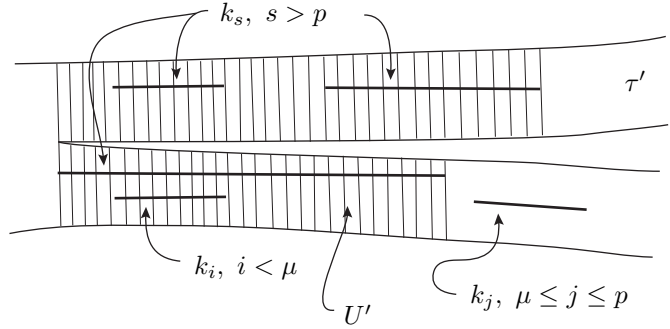}
	\end{center}
	\caption{Splitting to straighten shortcuts}
	\label{splittingshortcuts}
	\end{figure}

	Now let $U'$ denote that part of $U_0$ contained in $\tau'$, and note that the interior of $U'$ does not contain any switches of $\tau'$.

	Re-index the original collection of shortcuts $k_1,\dots,k_p$ so that those contained in $U_0$ are exactly the $k_i$ with $i<\mu$, and note that the remaining $k_{\mu},\dots,k_p$ are $\epsilon$-shortcuts for $\phi(\tau')$. We then complete these by a family of shortcuts $k_{p+1},\dots,k_{\nu}$ for $\phi(\tau')$
contained in $U'$ chosen so that the union of the $R'(k_j)$ is maximal with respect to the following properties:
	\begin{enumerate}
	\item the $R'(k_j)$ have disjoint interiors (for $j>p$, and thus for $j\geq\mu$),
	\item for all $i<\mu$, each component of $R_0(k_i)\cap U'$ is contained in some $R'(k_j)$ with $j>p$.
	\end{enumerate}

	We claim that the family $\{k_{\mu},\dots,k_{\nu}\}$ has the desired property, namely, that any other $\epsilon$-shortcut for $\phi(\tau')$ has at least one endpoint in the interior of some $R'(k_i)$ for $\mu\leq i\leq\nu$. To see this, suppose that $k$ is an $\epsilon$-shortcut for $\phi(\tau')$. Then $k$ is also an $\epsilon$-shortcut for $\phi(\tau_0)$. Suppose $\partial k$ avoids the interiors of the $R'(k_i)$ for $\mu\leq i\leq p$, and note that for these $i$ we have $R'(k_i)=R_0(k_i)$. Thus by the maximality condition on the original collection $k_1,\dots,k_p$ of $\epsilon$-shortcuts for $\tau_0$, either some endpoint of $k$ is contained in the interior of some $R_0(k_i)$ with $i<\mu$ or $k$ crosses a switch of $\tau_0$, and is therefore contained in $U_0\cap\tau'=U'$. In either case, we deduce that some endpoint of $k$ is contained in some $R'(k_j)$ with $j>p$ (by property (2) in the first case and by the maximality condition in the second case).
	\end{proof}

Here is the lemma that allows for the initial splitting used in the proof above.
\begin{lem}
Suppose $\tau$ is a weighted train track, no splitting of which cuts off an annular component. For any map $\phi$ adapted to $\tau$ there is a splitting $\tau'$ of $\tau$ so that $\ell(\phi(b))>C$ for every branch $b$ of $\tau'$, where $C$ is arbitrary.
\end{lem}

\begin{proof}
We construct a splitting of $\tau$ whose minimal branch length is $3/2$ times that of $\tau$. The result then follows by repeating this procedure.

Note that any splitting $\tau'$ of $\tau$, every switch of which is a switch of $\tau$, has minimal branch length bounded below by that of $\tau$. This is because every branch of $\tau'$ runs between two switches of $\tau$, and so runs at least the length of some branch of $\tau$.

Say a branch $b$ of $\tau$ is \emph{thick} if both its vertical boundary components contain corners in their interiors. If there are any homotopically nontrivial closed trainpaths carried by only thick branches (i.e., if there are any loops of thick branches), we split the branches over which such a path runs, cutting along arcs from corners of thick branches to existing switches of thick branches. By splitting far enough, we may remove all loops of thick branches (this follows from the fact that we have already assumed all annular portions of $\tau$ to have been removed). By the comment above, this does not decrease the minimum branch length.

We now assume $\tau$ has no loops of thick branches. Next, split $\tau$ along arcs from each corner to the midpoint of the first thick branch encountered. This splitting has minimal branch length at least that of $\tau$. To see this, note first that any branch of $\tau$ that is not thick will increase in length by at least half the length of one of its adjacent branches. Second, suppose there is a connected strip of $m$ thick branches of $\tau$. After splitting as described, we obtain a connected strip of $m-1$ thick branches, one for each switch between two thick branches. Each of these new thick branches has length equal to the average of the two thick branches that were adjacent to the corresponding switch. In particular, this splitting procedure does not decrease the minimum branch length. Moreover, if the longest strip of connected thick branches of $\tau$ has length $m$, then repeating this splitting procedure $m$ times will produce a train track with minimal branch length at least $3/2$ that of $\tau$. The result follows.
\end{proof}

\subsection{Proof of Proposition~\ref{prop:ending}}\label{ending}

We are given a homotopy equivalence $\phi\co S\to M$ from a closed hyperbolic surface $S$ into $M$, and a unit-length geodesic measured lamination $\alpha$ on $S$. We will prove Proposition~\ref{prop:ending} by assuming that there is some $\kappa>0$ so that $\ell(\phi'(\alpha))\geq\kappa$ for all $\phi'$ homotopic to $\phi$, and showing that conclusion (1) of the Proposition must hold.

We begin with a train track $\tau_0$ on $S$ fully carrying $\alpha$ with the property that any compact leaves of $\alpha$ are carried by annular components of $\tau_0$. We then homotope $\phi$ to $\phi_0$ adapted to $\tau_0$ carrying $\alpha$ so that these annuli are mapped to minimal edge-loops in $M$. Now we inductively define a sequence of splittings $\tau_n$ of $\tau_0$ and a sequence of maps $\phi_n\co S\to M$ adapted to $\tau_n$, so that $\phi_n$ and $\tau_n$ are obtained from $\phi_{n-1}$ and $\tau_{n-1}$ by the procedure of straightening shortcuts described above. Note that we may do this while ensuring that the maps all agree on the compact leaves of $\alpha$. In particular, if every component of $\alpha$ is a closed curve, then for all $n$ we may set $\tau_n=\tau_0$ and $\phi_n=\phi_0$.

Recall that we may cover the projective tangent bundle of $S$ with a finite number of compact flow boxes $B_i$, $1\leq i\leq q$ with disjoint interiors so that for $1\leq i\leq p$, all arcs of $B_i\cap\mathscr{F}$ project into $S$ as arcs transverse to ties in $\tau_{n}$, while for $p<i\leq q$, arcs of $B_i\cap\mathscr{F}$ are disjoint from $\alpha$ in $PT(S)$ (see Section~\ref{currents}).

	As usual, for any closed geodesic $\gamma\subset S$, we let $\gamma^*$ denote a minimal edge-loop in $M$ freely homotopic in $M$ to $\phi_n(\gamma)$. Note that because we are not assuming $\gamma$ to be simple, we may not ensure that $\gamma$ is entirely contained in $\tau_n$ by simply assuming $\gamma/\ell_S(\gamma)$ sufficiently close to $\alpha$. In other words, we may not assume that the lift of $\gamma$ to $PT(S)$ is disjoint from the boxes $B_{p+1},\dots,B_q$, and so there is generally some portion of $\gamma$ over whose image under $\phi_n$ we have no control. Let $\gamma\cap B$ denote that portion of $\gamma$ which lifts to $B$ in $PT(S)$. Given a closed geodesic $\gamma$ on $S$, we define $\gamma_n$ to be an edge-loop in $M^{(1)}$ obtained by replacing the image under $\phi_n$ of the components of $\gamma\cap B_i$ with $i>p$ with a minimal edge-path in $M$ to which it is homotopic relative to its endpoints.

\begin{figure}[ht]
	\begin{center}
		\psfrag{0}{$\gamma_n^0$}
		\psfrag{1}{$\gamma_n^1$}
		\psfrag{2}{$\gamma_n^2$}
		\psfrag{3}{$\gamma_n^3$}
		\psfrag{4}{$\gamma_n^4$}
		\psfrag{n}{$\gamma_n$}
		\psfrag{g}{$\gamma^*$}
		\includegraphics{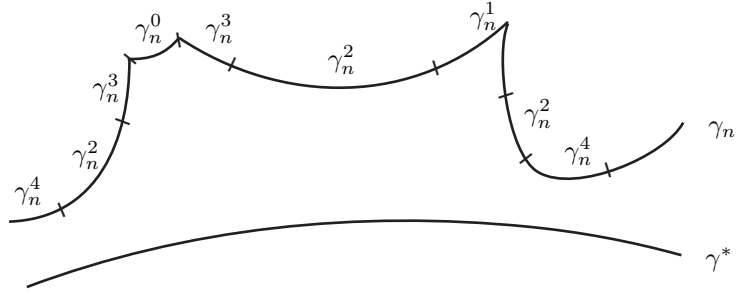}
	\end{center}
	\caption{The parts of $\gamma_n$}
	\label{gammas}
	\end{figure}

	We will partition the 1--cells of $\gamma_n$ into five classes, where the particular class in which a 1--cell sits is determined by the position of its endpoints relative to the other vertices of $\gamma^n$ as well as those of $\gamma_n^*$ (see Figure~\ref{gammas}). For this we need to fix a constant $\epsilon\in\Z$. We will specify below the conditions that $\epsilon$ must satisfy, but for now it suffices to suppose $\epsilon>1$.

	Define $\gamma_n^0$ to be those 1--cells of $\gamma_n$ either of whose endpoints is contained in \[\phi_n\big(\gamma\cap\bigcup_{i=p+1}^q B_i\big).\] Thus $\gamma_n^0$ is essentially that part of $\gamma_n$ which differs from $\phi_n(\gamma)$, namely those parts which are pulled tight in creating $\gamma_n$.

	Define $\gamma_n^4$ to be those 1--cells of $\gamma_n-\gamma_n^0$ which are within $\epsilon$ (in $M^{(1)}$) of some vertex of $\gamma^*$.

	To define $\gamma_n^1$, consider the splitting $\tau_{n+1}$ contained in $\tau_n$ obtained in the course of applying Lemma~\ref{lem:shortcuts} to straighten the shortcuts of $\phi_n(\tau_n)$. We define $\gamma_n^1$ to be all 1--cells $e$ of $\gamma_n-\gamma_n^0-\gamma_n^4$ for which there is some subarc $k$ of $\gamma$ with the following properties:
	\begin{enumerate}
	\item $k$ is an $\epsilon$-shortcut for $\phi_n(\tau_{n+1})$;
	\item $\phi_n(k)\subset\gamma_n-\gamma_n^0$;
	\item $e\cap \phi_n(k)\neq\emptyset$.
	\end{enumerate}
The third condition implies that $e$ and $\phi_n(k)$ share an endpiont.

	Define $\gamma_n^3$ to be those 1--cells in $\gamma_n-\gamma_n^0-\gamma_n^1-\gamma_n^4$ from which there is an edge-path $\lambda$ to some edge in $\gamma_n^0$ so that $\ell(\lambda)\leq\epsilon$ and $\lambda$ is homotopic with fixed endpoints to a subarc of $\gamma_n$. Essentially $\gamma_n^3$ consists of those edges from which there is, in some sense, an $\epsilon$-shortcut into $\gamma_n^0$.

	Finally, define $\gamma_n^2$ to be all remaining unassigned edges of $\gamma_n$. Thus edges in $\gamma_n^2$ are essentially those vertices of $\gamma_n-\gamma_n^0$ which are greater than a distance $\epsilon$ from $\gamma^*$, and from which there is no $\epsilon$-shortcut (carried by $\phi_n(\tau_{n+1})$) into $\gamma_n^1$ or $\gamma_n^0$.

	We aim to show that for $n$ large enough and $\overline{\gamma}$ close enough to $\alpha$, each of $\gamma_n^0$, $\gamma_n^1$, $\gamma_n^2$, and $\gamma_n^3$ is short, relative to the length of $\phi_n(\gamma)$, forcing $\gamma_n^4$ to be most of $\phi_n(\gamma)$ . In the following three lemmas, we obtain estimates on the absolute lengths of $\gamma_n^1$, $\gamma_n-\gamma_n^0$, and $\gamma_n^3$.

\begin{lem}
There is a constant $V_{\epsilon}$, depending on $G$ and $\epsilon$, so that
$$
\ell(\gamma_n^1)\leq 2V_{\epsilon}\Big(\ell(\phi_n(\gamma))-\ell(\phi_{n+1}(\gamma))+\epsilon\sum_{b\subset\tau_n'}\gamma(b)\Big).
$$
\end{lem}

\begin{proof}
	In order to control the length of $\gamma_n^1$, we consider the relationship between its length and the decrease in length to be realized when straightening the shortcuts of $\phi_n(\tau_n)$ to obtain $\phi_{n+1}(\tau_{n+1})$. When applying Lemma~\ref{lem:shortcuts} to straighten the shortcuts of $\phi_n(\tau_n)$, we obtain a finite maximal family of $\epsilon$-shortcuts $\{k_{\mu},\dots,k_{\nu}\}$ with the property that any other $\epsilon$-shortcut for $\phi_n(\tau_{n+1})$ has at least one endpoint contained in some $R'(k_i)$. Let $\mathcal{I}$ denote the vertex ties of $\phi_n(\tau_{n+1})$ contained in the union of the $R'(k_i)$, and let $V_{\epsilon}$ denote the volume of an $\epsilon$-ball in $X^{(1)}$. (Note that $V_{\epsilon}<\infty$ because $G$ is finitely generated.)

	We want to bound the length of $\gamma_n^1$ in terms of the size of $\mathcal{I}$ and $V_{\epsilon}$. We do this by associating each vertex in $\gamma_n^1$ to a tie of $\mathcal{I}$ and bounding the number of such vertices associated to any particular such tie in terms of $V_{\epsilon}$. To this end, fix a tie $t\in\mathcal{I}$. There is a collection of subarcs $\kappa_1,\dots,\kappa_s$ of $\gamma_n$ beginning at $t$ and forming $\epsilon$-shortcuts for $\phi_n(\tau_{n+1})$, as described in the definition of $\gamma_n^1$. We will call such arcs \emph{$(\gamma,t)$-spurs}. For each $i$, let $t_i$ denote the other endpoint of $\kappa_i$. We have no uniform control over the number of such $\kappa_i$ issuing from $t$, and so no control over the number of $t_i$ associated with $t$ in this way. We can, however, control the number of such $t_i$ which need to be counted as associated with $t$ because they will not be counted as associated to any other $t'\in\mathcal{I}$.

	Suppose $\kappa_i$ is a $(\gamma,t)$-spur having the property that if $\kappa_j$ is any other $(\gamma,t)$-spur issuing from $t$ in the same direction as $\kappa_i$ with $\phi_n(t_i)=\phi_n(t_j)$, then $\kappa_i\subset\kappa_j$. Then $\kappa_j-\kappa_i$ is a subarc of $\gamma$ carried by $\tau_{n+1}$ and with $\phi_n(\kappa_j-\kappa_i)$ a non-minimal edge-path in $M$ joining $\phi_n(t_i)$ to itself. In particular, $\kappa_j-\kappa_i$ is an $\epsilon$-shortcut for $\phi_n(\tau_{n+1})$, and thus either $t_i$ or $t_j$ is in $\mathcal{I}$.

	If $t_i\in\mathcal{I}$, then we need not count $t_j$ as associated with $t$, as it will be counted as associated to $t_i$ (or perhaps some other endpoint $t'$ of some other $(\gamma,t)$-spur $\kappa'$ with $\kappa_i\subset \kappa'\subset\kappa_j$). If $t_j\in\mathcal{I}$, then we clearly need not count it as associated with $t$. Thus if we associate to each $t\in\mathcal{I}$ the shortest $(\gamma,t)$-spur in either direction whose endpoint is mapped to a particular vertex in the $\epsilon$-ball about $\phi_n(t)$, we will have counted all points of $\gamma_n$ involved in computing the length of $\gamma_n^1$. We conclude that \[\ell(\gamma_n^1)\leq 2V_{\epsilon}\sum_{t\in\mathcal{I}}\gamma(t),\] where $\gamma(t)$ denotes the number of times that $\gamma$ intersects $t$ on $S$.

	Now we clearly have \[\ell(\phi_n(\gamma))=\sum_{t\notin\mathcal{I}}\gamma(t)+
\sum_{t\in\mathcal{I}}\gamma(t),\] and because the shortcuts $k_{\mu},\dots,k_{\nu}$ used to define $\tau'$ are straightened in defining $\phi_{n+1}$, we have that \[\ell(\phi_{n+1}(\gamma))\leq\sum_{t\notin\mathcal{I}}\gamma(t)+\epsilon\sum_{\mu\leq i\leq\nu}\gamma(k_i),\] where $\gamma(k_i)$ is the number of times $\gamma$ intersects a tie of $R'(k_i)$. Subtracting, we obtain \[\ell(\phi_n(\gamma))-\ell(\phi_{n+1}(\gamma))\geq
\sum_{t\in\mathcal{I}}\gamma(t)-\epsilon\sum_{\mu\leq i\leq\nu}\gamma(k_i).\]

	Combining this inequality with the one above, we obtain \[\ell(\gamma_n^1)\leq 2V_{\epsilon}\Big(\ell(\phi_n(\gamma))-\ell(\phi_{n+1}(\gamma))+\epsilon\sum_{\mu\leq i\leq\nu}\gamma(k_i)\Big).\]

	Now because $\phi_n$ is tightly adapted to $\tau_n$, and because each $\epsilon$-shortcut $k_i$ for $\phi_n(\tau_{n+1})$ is also an $\epsilon$-shortcut for $\phi_n(\tau_n)$, we see that each $R'(k_i)$ must cross a switch of $\tau_n$. Thus because the $R'(k_i)$ have disjoint interiors, we have that \[\sum_{\mu\leq i\leq\nu}\gamma(k_i)\leq\sum_{b\subset\tau_n'}\gamma(b),\] where $\tau_n'$ is the maximal subset of $\tau$ containing no annular components. (Recall that there are no shortcuts if $\tau_n'=\emptyset$.) We now have \[\ell(\gamma_n^1)\leq 2V_{\epsilon}\Big(\ell(\phi_n(\gamma))-\ell(\phi_{n+1}(\gamma))+\epsilon\sum_{b\subset\tau_n'}\gamma(b)\Big).\]
\end{proof}

\begin{lem}
\[\ell(\gamma_n-\gamma_n^0)\geq\ell(\phi_n(\gamma))-c_n\sum_{p<j\leq q}\gamma(B_j),\]
\end{lem}

\begin{proof}
	Let $c_n$ majorize the length of the images under $\phi_n$ of the arcs of $B_j\cap\mathscr{F}$ with $j>p$. The result is immediate (recall that $\gamma(B)$ is exactly the number of subarcs of $\gamma$ which lift to leaves of $B\cap\mathscr{F}$).
\end{proof}

\begin{lem}
\[\ell(\gamma_n^3)\leq d_n\ell(\gamma_n^0)\leq c_nd_n\sum_{p<j\leq q}\gamma(B_j).\]
\end{lem}

\begin{proof}
	Note that for each vertex $v$ in $\gamma_n^3$ there is a corresponding vertex $v'$ in $\gamma_n^0$. There is a constant $d_n$ depending on $\phi_n$ so that for any geodesic arc $\beta$ on $S$ with $\phi_n(\beta)$ homotopic with fixed endpoints to an edge-path with length no more than $\epsilon$ we must have $\ell(\beta)\leq d_n$. The result follows.
\end{proof}

	Before considering $\gamma_n^2$, we first show that for $n$ sufficiently large and $\overline{\gamma}$ sufficiently close to $\alpha$, we may assume that $\gamma_n^2$ and $\gamma_n^4$ together account for an arbitrarily large percentage of the length of $\gamma_n$.

\begin{lem}
For any $t<1$, there is an $n$ sufficiently large so that
$$
\frac{\ell(\gamma_n^4)+\ell(\gamma_n^2)}{\ell(\phi_n(\gamma))}>t
$$
for any $\overline{\gamma}$ sufficiently close to $\alpha$.
\end{lem}

\begin{proof}
	We begin with the fact that \[\ell(\gamma_n^4)+\ell(\gamma_n^2)=\ell(\gamma_n-\gamma_n^0)-\ell(\gamma_n^1)-\ell(\gamma_n^3).\] Using the above estimates for the lengths of $\gamma_n^0$, $\gamma_n^1$, and $\gamma_n^3$, we obtain \[\frac{\ell(\gamma_n^4)+\ell(\gamma_n^2)}{\ell(\phi_n(\gamma))}\geq 1-c_n(d_n+1)\sum_{p<j\leq
q}\frac{\gamma(B_j)}{\ell(\phi_n(\gamma))}\] \[-2V_{\epsilon}\Big(1-\frac{\ell(\phi_{n+1}(\gamma))}{\ell(\phi_n(\gamma))}+\epsilon\sum_{b\subset\tau_n'}\frac{\gamma(b)}{\ell(\phi_n(\gamma))}\Big).\] This inequality holds for all $n$ and for all $\gamma$. We now consider what happens to the right side of this inequality as $\overline{\gamma}\to\alpha$ with $n$ fixed.

	First note that \[\frac{\gamma(B_j)}{\ell(\phi_n(\gamma))}=\frac{\gamma(B_j)/\ell(\gamma)}{\ell(\phi_n(\gamma))/\ell(\gamma)}=\frac{\overline{\gamma}(B_j)}{\ell(\phi_n(\overline{\gamma}))}.\] Thus we have for fixed $n$ that \[\lim_{\overline{\gamma}\to\alpha}\frac{\gamma(B_j)}{\ell(\phi_n(\gamma))}=\frac{\alpha(B_j)}{\ell(\phi_n(\alpha))}.\] Similarly we find that \[\lim_{\overline{\gamma}\to\alpha}\frac{\ell(\phi_{n+1}(\gamma))}{\ell(\phi_n(\gamma))}=\frac{\ell(\phi_{n+1}(\alpha))}{\ell(\phi_n(\alpha))}\] and \[\lim_{\overline{\gamma}\to\alpha}\frac{\gamma(b)}{\ell(\phi_n(\gamma))}=\frac{\alpha(b)}{\ell(\phi_n(\alpha))}\] for each branch $b$ of $\tau_n'$. Thus for any fixed $n$, we have that \[\lim_{\overline{\gamma}\to\alpha}\frac{\ell(\gamma_n^4)+\ell(\gamma_n^2)}{\ell(\phi_n(\gamma))}\geq 1-c_n(d_n+1)\sum_{p<j\leq q}\frac{\alpha(B_j)}{\ell(\phi_n(\alpha))}\]\[-2V_{\epsilon}\Big(1-\frac{\ell(\phi_{n+1}(\alpha))}{\ell(\phi_n(\alpha))}+\epsilon\sum_{b\subset\tau_n'}\frac{\alpha(b)}{\ell(\phi_n(\alpha))}\Big).\] By definition, $\alpha(B_j)=0$ for all $p<j\leq q$, while $\ell(\phi_n(\alpha))$ is bounded away from zero for all $n$. The inequality thus becomes \[\lim_{\overline{\gamma}\to\alpha}\frac{\ell(\gamma_n^4)+\ell(\gamma_n^2)}{\ell(\phi_n(\gamma))}\geq 1-2V_{\epsilon}\Big(1-\frac{\ell(\phi_{n+1}(\alpha))}{\ell(\phi_n(\alpha))}+\epsilon\sum_{b\subset\tau_n'}\frac{\alpha(b)}{\ell(\phi_n(\alpha))}\Big).\] Moreover this inequality holds for all $n$.

	Now set $\alpha'$ to be that portion of $\alpha$ carried by $\tau_n'$ (i.e., $\alpha'$ is the union of all components of $\alpha$ which are not simple closed curves). If $\ell(\alpha')=0$, then because $\ell(\phi_n(\alpha))$ is assumed bounded away from zero, we have \[\sum_{b\subset\tau_n'}\frac{\alpha(b)}{\ell(\phi_n(\alpha))}=0.\] On the other hand, suppose $\ell(\alpha')>0$, and set $m_n$ to be the minimum branch length of $\tau_n'$, so that for all $n$ we have \[0<m_n\sum_{b\subset\tau_n'}\alpha(b)\leq\sum_{b\subset\tau_n'}\ell(b)\alpha(b)=\ell(\alpha').\]  Note that because $\tau_n'$ carries no closed curves of $\alpha$, $m_n$ grows arbitrarily large with $n$. It follows that \[\lim_{n\to\infty}\sum_{b\subset\tau_n'}\alpha(b)=0.\] Finally note that for large $n$ we have \[\lim_{n\to\infty}\frac{\ell(\phi_{n+1}(\alpha))}{\ell_n(\phi_n(\alpha))}=1.\] The result follows.
\end{proof}

We next choose $t$.

	Let $L=\max(12\delta+2\epsilon,4\epsilon)$, and choose $t$ as follows. Suppose $A$ is a finite circularly ordered set and $B$ is a subset of $A$. We denote by $B_L$ the set of elements of $B$ which are contained in a connected sequence of length at least $L$. We choose $t<1$ large enough to ensure that for $|A|$ sufficiently large, we have that $|B|/|A|\geq t$ implies that $|B_L|/|A|>1/2$. (The choice of $1/2$ is somewhat arbitrary.) Note that we necessarily have $t>1/2$.

	Given this choice of $t$, we henceforth assume that $\overline{\gamma}$ is chosen close enough to $\alpha$ and $n$ is chosen large enough to ensure that \[\frac{\ell(\gamma_n^2)+\ell(\gamma_n^4)}{\ell(\phi_n(\gamma))}>t.\] We choose $t$ in this way because, rather than bound the relative length of all of $\gamma_n^2$, we will bound the relative length of those portions of $\gamma_n^2$ which may be connected to form segments with length at least $4\epsilon$.

\begin{figure}[ht]
	\begin{center}
		\psfrag{A}{$A\left\{\rule[-1.5cm]{0cm}{3cm}\right.$}
		\psfrag{B}{$A_{\epsilon/2}\left\{\rule[-.75cm]{0cm}{1.5cm}\right.$}
		\psfrag{C}{$C\subset\hat{\gamma}_n^2$}
		\psfrag{F}{$\overbrace{\rule{4cm}{0cm}}$}
		\psfrag{D}{$C_{\epsilon}'\subset\hat{\gamma}_{\epsilon/2}^2$}
		\psfrag{E}{$C'$}
		\psfrag{g}{$\gamma^*$}
		\psfrag{e}{$\gamma_{\epsilon/2}$}
		\psfrag{n}{$\gamma_n$}
		\includegraphics{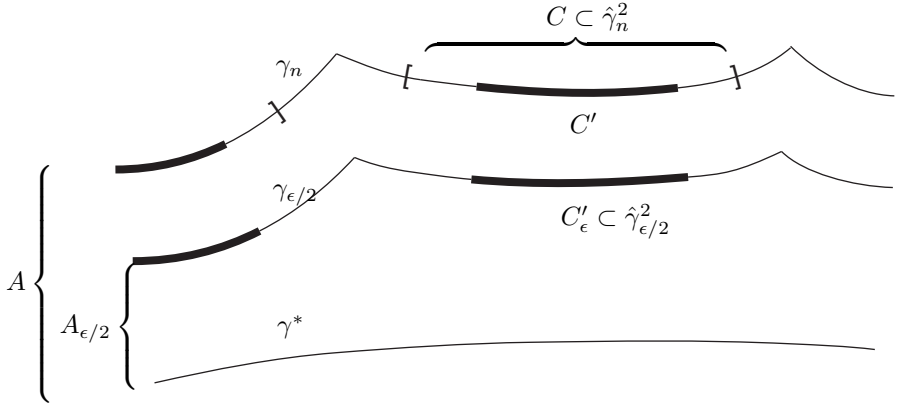}
	\end{center}
	\caption{Using area to bound the length of $\gamma_n^2$}
	\label{gamma2area}
	\end{figure}

	To begin, let $A_n\co S^1\times[0,1]\to M^{(2)}$ be a simplicial map of an annulus into $M^{(2)}$ so that $A(S^1\times\{0\})=\gamma^*$ and $A(S^1\times\{1\})=\gamma_n$. We pull back the metric from $M^{(2)}$ to $S^1\times[0,1]$ to produce a simplicial metric structure on the annulus, and note that the boundary components have lengths $\ell(\gamma^*)$ and $\ell(\gamma_n)$. For simplicity, we will refer to the annulus with this structure as $A_n$, and identify the boundary components of $A_n$ with $\gamma^*$ and $\gamma_n$ (see Figure~\ref{gamma2area}).

	Assume that $A_n$ is chosen to have minimal area. Let $A_{\epsilon/2}$ denote the complement in $A_n$ of the $\epsilon/2$-neighborhood (in $A_n$) of $\gamma_n$, and let $\gamma_{\epsilon/2}$ denote the boundary of this neighborhood in $A_n$. Let $\hat{\gamma}_n^2$ denote the union of edges in $\gamma_n$ both of whose endpoints are in $\gamma_n^2$. Let $\gamma_{\epsilon/2}^2$ denote those vertices $v$ in $\gamma_{\epsilon/2}$ with $d(v,\gamma_n)=d(v,\gamma_n^2)$, and similarly let $\hat{\gamma}_{\epsilon/2}^2$ denote the union of the edges of $\gamma_{\epsilon/2}$ both of whose endpoints are in $\gamma_{\epsilon/2}^2$ (see Figure~\ref{gamma2area}).

	Let $\mathcal{L}$ denote the union of the components of $\hat{\gamma}_n^2$ with length at least $L$, and let $C$ be a component of $\mathcal{L}$. Note that because $\gamma_n$ is piecewise geodesic, and any vertices around which $\gamma_n$ is not minimal must be contained in $\gamma_n^0$, $\gamma_n^1$, $\gamma_n^3$, or $\gamma_n^4$, it must be that $C$ corresponds to a minimal edge-path in $M^{(1)}$. Now consider the subsegment $C'$ of $C$ with endpoints a distance $\epsilon$ from those of $C$, and note then that $C'$ has length at least $\max(12\delta,2\epsilon)$. Corresponding to $C'$ is a subsegment $C_{\epsilon/2}'$ of $\hat{\gamma}_{\epsilon/2}^2$ consisting of those points in $\hat{\gamma}_{\epsilon/2}^2$ which are exactly a distance $\epsilon/2$ in $A_n$ from some point in $C'$. 

	Because the length of $C'$ is at least $12\delta$, it follows from the proof of Lemma~\ref{lem:expgrowth} that \[\ell(C_{\epsilon/2}')\geq c_12^{\epsilon/2}\,\ell(C').\] Moreover, if $\hat{C}$ is another component of $\hat{\gamma}_n^2$ with length at least $L$, then no point of $C_{\epsilon}$ lies within $\epsilon/2$ in $A_n$ of any point on the subsegment $\hat{C}'$ of $\hat{C}$ a distance $\epsilon$ from the endpoints of $\hat{C}$. For if $x\in C_{\epsilon}$ were such a point, lying within $\epsilon/2$ of $y\in C'$ and $z\in\hat{C}'$, then $d(y,z)\leq d(y,x)+d(x,z)\leq\epsilon$, while the subarc of $\gamma_n$ joining $y$ to $z$ has length at least $2\epsilon$. Thus $y$ and $z$ lie either in $\gamma_n^0$, $\gamma_n^1$, or $\gamma_n^4$, contradicting the choice of $C$ and $\hat{C}$.

	It follows from the discussion above, along with the fact that $\ell(C')\geq\ell(C)/2$, that \[\ell(\gamma_{\epsilon})\geq(c_1/2)\ell(\mathscr{L})2^{\epsilon/2}.\] From Lemma~\ref{lem:isoannulus} we find that \[\ell(\gamma_{\epsilon})\leq\text{area}(N_1(\gamma_{\epsilon}))\leq \text{area}(A)\leq K\big(\ell(\gamma_n)+\ell(\gamma^*)\big)\leq 2K\ell(\gamma_n),\] where $N_1(\gamma_{\epsilon})$ denotes the combinatorial 1--neighborhood of $\gamma_{\epsilon}$. Combining this with the inequality above, we find that \[\ell(\mathscr{L})\leq\frac{4K}{c_1}\,\ell(\gamma_n)2^{-\epsilon/2}.\]

	Now recall we are assuming $t$ fixed, $n$ large enough, and $\overline{\gamma}$ close enough to $\alpha$ so that \[\ell(\gamma_n^2)+\ell(\gamma_n^4)>t\ell(\phi_n(\gamma)).\] By our choice of $t$, we also have \[\ell(\mathscr{L})>\frac{1}{2}\,\ell(\gamma_n^2).\] Combining these with the inequality obtained above for $\ell(\mathscr{L})$, we find that
\[t\,\ell(\phi_n(\gamma))<\frac{8K}{c_1}\,\ell(\gamma_n)2^{-\epsilon/2}+\ell(\gamma_n^4).\]

	Note that as $\overline{\gamma}\to\alpha$, we have $\dfrac{\ell(\gamma_n)}{\ell(\phi(\gamma))}\to 1$. It follows that for $\overline{\gamma}$ sufficiently close to $\alpha$, we have \[\frac{\ell(\gamma_n^4)}{\ell(\phi_n(\gamma))}>t-\frac{8K}{c_1}2^{-\epsilon/2}.\] Thus by choosing $t$ sufficiently close to $1$ and $\epsilon$ sufficiently large, we can ensure that, for $n$ sufficiently large, the ratio $\ell(\gamma_n^4)/\ell(\phi_n(\gamma))$ is arbitrarily close to $1$ for $\overline{\gamma}$ sufficiently close to $\alpha$. This completes the proof of Proposition~\ref{prop:ending}.

\section{Proof of Theorem~\ref{thm:main}}\label{proof}

As outlined in the introduction, the first step in the proof of Theorem~\ref{thm:main} is the following lemma (cf.~\cite[Prop $2.3$]{bonahon:ends}).

	\begin{lem}\label{lem:curvesexit}
	Suppose $b$ is a geometrically infinite end of $M$. Then there is a sequence of closed curves $\alpha_i$ in $M$ so that $\alpha_i^*\to b$.
	\end{lem}

	\begin{proof}
	All distances will be measured in the 1--skeleton of the appropriate space.

	Let $U$ be a neighborhood of $b$ and choose a neighborhood $U_0$ of $b$ contained in $U$ and a distance at least $2c_2$ from the complement of $U$ in $M$, where $c_2$ is the constant from Lemma~\ref{lem:nbhd}. We may assume (cf. Section~\ref{basic}) that $U_0$ is connected and has compact boundary $\partial U_0\subset M^{(2)}$. To prove the lemma, it suffices to find a curve $\alpha$ in $M$ with $\alpha^*\subset U$.

	Because $b$ is geometrically infinite, there are infinitely many minimal edge-loops $\sigma_i^*$ intersecting $U_0$ nontrivially. If any one of these is contained in $U_0$, we are finished. If not, these geodesics give rise to an infinite family of locally minimal arcs $\{\overline{\sigma}_i\}$ contained in $U_0$ and with endpoints on vertices in $\partial U_0$ (see Figure~\ref{exitend}). Because $\partial U_0$ is compact, we may subsequence so that $\partial\overline{\sigma}_i=\partial\overline{\sigma}_j$ for all $i,j$. Because the $\overline{\sigma}_i$ intersect arbitrary neighborhoods of $b$, we may further assume that at least one point of $\overline{\sigma}_j$ is a distance at least $3\delta$ from all $\overline{\sigma}_i$ with $i<j$.

	\begin{figure}[ht]
	\begin{center}
		\psfrag{U}{$U$}
		\psfrag{u}{$U_0$}
		\psfrag{a}{$\alpha$}
		\includegraphics{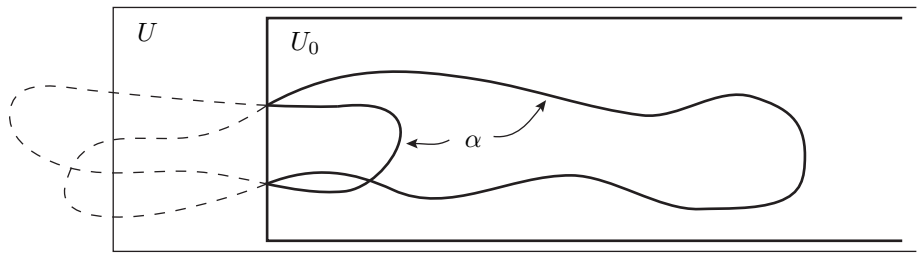}
	\end{center}
	\caption{Finding curves exiting an end}
	\label{exitend}
	\end{figure}

	Fix $i$ and $j$ and define the closed curve $\alpha$ to be the composition $\overline{\sigma}_i\circ\overline{\sigma}_j$. Note that $\alpha$ is not null-homotopic in $M$, since if it were, then $\overline{\sigma}_i$ and $\overline{\sigma}_j$ would lift to $X^{(1)}$ to two minimal edge-paths with the same endpoints, one of which contains a point a distance at least $3\delta$ from the other, contradicting the fact that geodesic triangles in $X^{(1)}$ are $\delta$-thin.

	Applying Lemma~\ref{lem:nbhd}, we see that any minimal edge-loop $\alpha^*$ corresponding to $\alpha$ must lie within a $2c_2$ neighborhood of $\alpha$, and thus lie within $U$, as required.
	\end{proof}

	We may choose the $\alpha_i$ to be geodesics on $S$ with respect to some arbitrary hyperbolic metric, and then subsequence so that $\alpha_i/\ell_S(\alpha_i)\to\alpha_{\infty}$ for some unit length geodesic current $\alpha_{\infty}$.

	\begin{lem}\label{lem:lam}
	Suppose $\alpha_i$ is a sequence of geodesics on $S$ with $\alpha_i/\ell_S(\alpha_i)$ converging in $\mathscr{C}(S)$ to a geodesic current $\alpha_{\infty}$. If the minimal edge-loops $\alpha_i^*$ corresponding to $\phi(\alpha_i)$ leave compact sets, then $\alpha_{\infty}$ is in fact a measured lamination on $S$.
	\end{lem}

	\begin{proof}
	By Lemma~\ref{lem:intnum} we have that if $D_i=d(\phi(\alpha_i),\alpha_i^*)$, then
\[\lim_{i\to\infty}\text{int}\Big(\frac{\alpha_i}{\ell(\phi(\alpha_i))},\frac{\alpha_i}{\ell(\phi(\alpha_i))}\Big)\leq\lim_{i\to\infty}C\,2^{-D_i}=0.\] Because $S$ is compact and $\phi$ is proper, the ratio $\ell(\phi(\alpha_i))/\ell_S(\alpha_i)$ is bounded above by a constant independent of $i$. Thus \[\lim_{i\to\infty}\frac{\text{int}(\alpha_i,\alpha_i)}{\ell_S(\alpha_i)^2}\to 0.\] By continuity of intersection number, it follows that \[\text{int}(\alpha_{\infty},\alpha_{\infty})=0.\] By Lemma~\ref{lem:lam}, $\alpha_{\infty}$ is a lamination.
	\end{proof}

Now because $\alpha_{\infty}$ is a lamination, there is a sequence of simple closed geodesics $\gamma_j$ on $S$ with $\gamma_j/\ell_S(\gamma_j)\to\alpha_{\infty}$ in $\mathscr{C}(S)$. As the $\alpha_i^*$ leave compact sets of $M$, we know that $\alpha_{\infty}$ does not satisfy conclusion (1) of Proposition~\ref{prop:ending}. It follows that as $\phi$ varies over its homotopy class, $\ell(\phi(\alpha_{\infty}))$ may be made arbitrarily small. We now show that because the $\alpha_i^*$ leave compact sets of $M$, the $\gamma_j^*$ at least are not contained in any compact set of $M$.

	\begin{lem}\label{lem:leave}
	With notation as above, suppose $\gamma_j$ are simple closed geodesics on $S$ with $\gamma_j/\ell_S(\gamma_j)\to\alpha_{\infty}$ for some unit length measured lamination $\alpha_{\infty}$ satisfying conclusion (2) of Proposition~\ref{prop:ending}. Then there is no compact set of $M$ containing all the minimal edge-loops $\gamma_j^*$ corresponding to the $\phi(\gamma_j)$.
	\end{lem}

	\begin{proof}
	This is the argument of~\cite[\S $6.3$]{bonahon:ends}.

	Suppose otherwise, so that the $\gamma_i^*$ are all entirely contained in some compact set $K$. Because $K$ is compact the ratio $\ell(\gamma_i^*)/\ell_S(\gamma_i)$ is bounded below by some constant $c$ depending only on $K$ and $\delta$. Thus the ratio $\ell(\gamma_i')/\ell_S(\gamma_i)$ is similarly bounded for all edge-loops $\gamma_i'$ homotopic to $\gamma_i$ in $M$.

	In particular, for any map $\phi'$ homotopic to $\phi$, we have that for all $i$, \[\ell(\phi'(\gamma_i))/\ell_S(\gamma_i)\geq c.\] By letting $i$ tend to $\infty$, we find that \[\ell(\phi'(\alpha))\geq c^{-1}>0.\] That this is true for any $\phi'$ contradicts the assumption that $\alpha_{\infty}$ satisfies conclusion (2) of Proposition~\ref{prop:ending}. Thus there is no compact set of $M$ containing $\gamma_i^*$ for all $i$.
	\end{proof}

\bibliography{bib}
\bibliographystyle{alpha}
\end{document}